# A SOLUTION TO THE BAER SPLITTING PROBLEM

LIDIA ANGELERI HÜGEL, SILVANA BAZZONI, AND DOLORS HERBERA

ABSTRACT. Let $R$ be a commutative domain. We prove that an $R$-module $B$ is projective if and only if $\mathrm{Ext}^1_R(B,T) = 0$ for any torsion module $T$. This answers in the affirmative a question raised by Kaplansky in [16].

A module $B$ over a commutative domain $R$ is called a Baer module in case $\mathrm{Ext}^1_R(B,T) = 0$ for every torsion $R$-module $T$. This definition goes back to 1936 when R. Baer [3] posed the question of characterizing the class of all abelian (torsion-free) groups $G$ such that any extension of $G$ with a torsion group splits. In the language of homological algebra, the problem asks which groups $G$ satisfy $\mathrm{Ext}^1_{\mathbb{Z}}(G,T) = 0$ for all torsion groups $T$. Baer proved that every countably generated group $G$ with this property must be free [3, Theorem 8.6 and Footnote 11 p. 781].

In 1961 Rotman [19] introduced the terminology of Baer groups or B-groups and put this problem, together with the Whitehead problem, in the more general setting of describing, for a given class of abelian groups $\mathcal{S}$, the groups $B$ satisfying that $\mathrm{Ext}^1_{\mathbb{Z}}(B,S) = 0$ for any abelian group $S \in \mathcal{S}$.

In 1962 Kaplansky [16] considered the case of modules over commutative domains. He raised the question whether Baer modules are projective. Using what now are well known tools of homological algebra, he proved that Baer modules are flat, hence torsion-free, modules of projective dimension at most one.

The answer to the original problem raised by Baer was only given in 1969, when Griffith [10] proved that the only Baer groups are the free groups. Grimaldi [12] later generalized the result to modules over Dedekind domains proving that the only Baer modules over such domains are the projective modules.

A real breakthrough in the study of the structure of Ext was made by Shelah in [20], showing that set theoretic methods are essential in this area. Following this track, in 1988 Eklof and Fuchs [5] used a version of Shelah's Singular Compactness Theorem to prove that Baer modules over valuations domains are free. In [8] Eklof, Fuchs and Shelah, generalized

*Date*: January 18, 2006.
*Key words and phrases.* Baer modules, Mittag-Leffler inverse sytems
2000 *Mathematics Subject Classification*. Primary: 13C05; 16E30. Secondary: 13G05;16D40.
First and second author supported by Università di Padova (Progetto di Ateneo CDPA048343 "Decomposition and tilting theory in modules, derived and cluster categories"). First and third author supported by the DGI and the European Regional Development Fund, jointly, through Project MTM2005-00934. Third author supported by the Comissionat per Universitats i Recerca of the Generalitat de Catalunya and by the "Gruppo Nazionale per le Strutture Algebriche, Geometriche e le loro Applicazioni" (Italy). Part of this paper was written while the third author was visiting the universities in Padova and in Varese, she wants to thank her hosts for their hospitality.

the tools used in [5] to arbitrary domains, and they proved a crucial reduction theorem. Namely, they showed that a module $B$ over an arbitrary domain is a Baer module if and only if it is the union of a well ordered continuous ascending chain $(B_\alpha \mid \alpha < \lambda)$ of submodules such that the factors $B_{\alpha+1}/B_\alpha$ are countably generated Baer modules. This reduces the problem of showing that Baer modules are projective to the countably generated case. For these results, as well as for a general account on the problem of studying the structure of Ext, we refer to [9, Chapter XVI §8] and [6, Chapter XII §3].

Since then, to our knowledge, the only substantial progress concerning the Baer problem was made by Griffith [11], who showed that Baer modules over local noetherian regular domains are free.

In the present paper, we show that every countably generated Baer module over an arbitrary commutative domain is projective. Hence, by the result of Eklof, Fuchs and Shelah, it follows that all Baer modules are projective. This solves the general problem raised by Kaplansky, see also [6, Open Problem F.2] and [9, Problem 60].

Let us sketch the idea of the proof. Let $B$ be a countably generated Baer module over a commutative domain $R$. Then it is well known that $B$ is countably presented and flat [16]. So, there is a countable direct system

$$F_1 \xrightarrow{f_1} F_2 \xrightarrow{f_2} F_3 \to \ldots \to F_n \xrightarrow{f_n} F_{n+1} \to \ldots$$

of finitely generated free modules $F_n$ such that $B = \varinjlim F_n$, see [14].

As a first step, we use recent work of the second and third author [4] to translate the vanishing of Ext into a Mittag-Leffler condition on inverse systems. More precisely, we show that $B$ is a Baer module iff for any torsion module $M$ the functor $\mathrm{Hom}(-, M)$ maps the direct system $F_1 \xrightarrow{f_1} F_2 \xrightarrow{f_2} F_3 \to \ldots$ into a Mittag-Leffler tower $(\mathrm{Hom}_R(F_n, M), \mathrm{Hom}_R(f_n, M))_{n \in \mathbb{N}}$, that is, an inverse system

$$\ldots \to \mathrm{Hom}(F_3, M) \xrightarrow{\mathrm{Hom}(f_2, M)} \mathrm{Hom}(F_2, M) \xrightarrow{\mathrm{Hom}(f_1, M)} \mathrm{Hom}(F_1, M)$$

satisfying the Mittag-Leffler condition of [13]. Similarly, we see that $B$ is a projective module iff the same holds true for the functor $\mathrm{Hom}(-, R)$.

Now the core of the proof consists in studying closure properties of the class of all modules $M$ that turn the direct system $F_1 \xrightarrow{f_1} F_2 \xrightarrow{f_2} F_3 \to \ldots$ into a Mittag-Leffler tower $(\mathrm{Hom}_R(F_n, M), \mathrm{Hom}_R(f_n, M))_{n \in \mathbb{N}}$. This will allow us to deduce that also the regular module $R$ belongs to this class, which completes the proof.

The paper is organized as follows. Section 1 is devoted to the general notion of a Mittag-Leffler tower. In Section 2, we consider a countable direct system of the form $C_1 \xrightarrow{f_1} C_2 \xrightarrow{f_2} C_3 \to \ldots$ together with a module $M$ over an arbitrary ring, and we give criteria for $(\mathrm{Hom}_R(C_n, M), \mathrm{Hom}_R(f_n, M))_{n \in \mathbb{N}}$ being a Mittag-Leffler tower. Finally, in Section 3, we apply our investigations to countably generated Baer modules over commutative domains and prove our main result.

Our rings are associative, have an identity, and they are not necessarily commutative unless stated otherwise. Modules are unital.

We wish to thank the referee for suggesting the use of Lemma 2.7 in order to simplify the proof of Lemma 2.8.

## 1. MITTAG-LEFFLER TOWERS

For a given set $I$, let $\{M_i\}_{i\in I}$ and $\{N_i\}_{i\in I}$ be two families of right modules over a ring $R$, and let $\{\gamma_i\colon M_i \to N_i\}_{i\in I}$ be a family of module homomorphisms. Then there is a module homomorphisms $\Gamma\colon \prod_{i\in I} M_i \to \prod_{i\in I} N_i$ defined by $\Gamma((m_i)_{i\in I}) = (\gamma_i(m_i))_{i\in I}$ for any $(m_i)_{i\in I} \in \prod_{i\in I} M_i$. Moreover, $\Gamma$ induces by restriction a homomorphism $\Gamma'\colon \bigoplus_{i\in I} M_i \to \bigoplus_{i\in I} N_i$. We will call $\Gamma$ and $\Gamma'$ the *diagonal maps induced by* $\{\gamma_i\}_{i\in I}$, or simply diagonal maps.

It is immediate to check that the composition of diagonal maps is a diagonal map. More precisely,

**Lemma 1.1.** *Let $R$ be a ring and $I$ a set. Let $\{K_i\}_{i\in I}$, $\{M_i\}_{i\in I}$ and $\{N_i\}_{i\in I}$ be families of right $R$-modules, and let $\{\lambda_i\colon K_i \to M_i\}_{i\in I}$ and $\{\gamma_i\colon M_i \to N_i\}_{i\in I}$ be families of maps. Denote by $\Lambda\colon \prod_{i\in I} K_i \to \prod_{i\in I} M_i$, $\Lambda'\colon \bigoplus_{i\in I} K_i \to \bigoplus_{i\in I} M_i$, $\Gamma\colon \prod_{i\in I} M_i \to \prod_{i\in I} N_i$ and $\Gamma'\colon \bigoplus_{i\in I} M_i \to \bigoplus_{i\in I} N_i$ the induced diagonal maps.*

*Then $\Gamma\Lambda\colon \prod_{i\in I} K_i \to \prod_{i\in I} N_i$ and $\Gamma'\Lambda'\colon \bigoplus_{i\in I} K_i \to \bigoplus_{i\in I} N_i$ are the diagonal maps induced by $\{\gamma_i\lambda_i\}_{i\in I}$.*

*Moreover, the following statements are equivalent*

(1) $\Gamma\Lambda(\prod_{i\in I} K_i) = \Gamma(\prod_{i\in I} M_i)$,
(2) *for any $i \in I$, $\gamma_i\lambda_i(K_i) = \gamma_i(M_i)$,*
(3) $\Gamma'\Lambda'(\bigoplus_{i\in I} K_i) = \Gamma'(\bigoplus_{i\in I} M_i)$.

PROOF. All the statements follow easily from the fact that $\Lambda$, $\Gamma$ and their restrictions, $\Lambda'$ and $\Gamma'$, respectively, are defined componentwise. ∎

A (countable) *tower* $\mathcal{T}$ of right $R$-modules consists of a sequence of modules $(H_n)_{n\in\mathbb{N}}$ and a sequence of morphisms

$$\ldots \to H_{n+1} \xrightarrow{\lambda_n} H_n \to \ldots \to H_3 \xrightarrow{\lambda_2} H_2 \xrightarrow{\lambda_1} H_1$$

We will use the notation $\mathcal{T} = (H_n, \lambda_n)_{n\in\mathbb{N}}$. Note that a tower of right $R$-modules is a representation in Mod-$R$ of the quiver

$$A_\infty: \quad \cdots \to \bullet \to \bullet \cdots \to \bullet \to \bullet$$

If $\mathcal{T}_i = (H_n^i, \lambda_n^i)_{n\in\mathbb{N}}$, $i = 1, 2$, are towers, then by a morphism of towers $f\colon \mathcal{T}_1 \to \mathcal{T}_2$ we mean a sequence of module homomorphisms $(f_n\colon H_n^1 \to H_n^2)_{n\in\mathbb{N}}$ satisfying $f_n\lambda_n^1 = \lambda_n^2 f_{n+1}$, for any $n \in \mathbb{N}$. If each $f_n$ is an isomorphism then we say that $\mathcal{T}_1$ and $\mathcal{T}_2$ are isomorphic towers.

The notions of direct sum and product of towers are crucial for our investigations. For a given set $I$, let $\{\mathcal{T}_i = (H_n^i, \lambda_n^i)_{n\in\mathbb{N}}\}_{i\in I}$ be a family of towers. Then the product and the direct sum of $\{\mathcal{T}_i\}_{i\in I}$ are

$$\prod_{i\in I} \mathcal{T}_i = (\prod_{i\in I} H_n^i, \Lambda_n)_{n\in\mathbb{N}} \text{ and } \bigoplus_{i\in I} \mathcal{T}_i = (\bigoplus_{i\in I} H_n^i, \Lambda'_n)_{n\in\mathbb{N}},$$

where, for each $n \in \mathbb{N}$,

$$\Lambda_n\colon \prod_{i\in I} H_{n+1}^i \to \prod_{i\in I} H_n^i \text{ and } \Lambda'_n\colon \bigoplus_{i\in I} H_{n+1}^i \to \bigoplus_{i\in I} H_n^i$$

are the diagonal maps induced by $\{\lambda_n^i\}_{i \in I}$, respectively.

If $\mathcal{T}_i = \mathcal{T}_j$, for any $i, j \in I$, we write $\prod_{i \in I} \mathcal{T}_i = \mathcal{T}^I$ and $\bigoplus_{i \in I} \mathcal{T}_i = \mathcal{T}^{(I)}$.

We recall the definition of a *Mittag-Leffler tower*, see [13] or [21, Definition 3.5.6].

**Definition.** A tower of right $R$-modules, $\mathcal{T} = (H_n, \lambda_n)_{n \in \mathbb{N}}$, satisfies the *Mittag-Leffler condition* if, for every $m \in \mathbb{N}$, the chain of submodules of $H_m$

$$H_m \supseteq \lambda_m(H_{m+1}) \supseteq \cdots \supseteq \lambda_m \lambda_{m+1} \cdots \lambda_{m+n-1}(H_{m+n}) \supseteq \cdots$$

is stationary, that is, if for each $m \in \mathbb{N}$, there exists $l(m) > m$ such that

$$\lambda_m \cdots \lambda_k(H_{k+1}) = \lambda_m \cdots \lambda_{l(m)-1}(H_{l(m)})$$

for any $k \geq l(m)$. If we are interested in recording the values $l(m)$, we say that the tower $\mathcal{T}$ satisfies the Mittag-Leffler condition with respect to the sequence $(l(m))_{m \in \mathbb{N}}$.

**Proposition 1.2.** *Let $R$ be a ring and $I$ a set. Let $\{\mathcal{T}_i = (H_n^i, \lambda_n^i)_{n \in \mathbb{N}}\}_{i \in I}$ be a family of towers. Let $(l(m))_{m \in \mathbb{N}}$ be a sequence of integers such that $l(m) > m$.*

*Then the following statements are equivalent.*

(1) $\prod_{i \in I} \mathcal{T}_i$ *satisfies the Mittag-Leffler condition with respect to the sequence* $(l(m))_{m \in \mathbb{N}}$,
(2) *every $\mathcal{T}_i$ satisfies the Mittag-Leffler condition with respect to the sequence* $(l(m))_{m \in \mathbb{N}}$,
(3) $\bigoplus_{i \in I} \mathcal{T}_i$ *satisfies the Mittag-Leffler condition with respect to the sequence* $(l(m))_{m \in \mathbb{N}}$.

PROOF. Set $\prod_{i \in I} \mathcal{T}_i = (\prod_{i \in I} H_n^i, \Lambda_n)_{n \in \mathbb{N}}$ and $\bigoplus_{i \in I} \mathcal{T}_i = (\bigoplus_{i \in I} H_n^i, \Lambda_n')_{n \in \mathbb{N}}$. Then $\prod_{i \in I} \mathcal{T}_i$ satisfies the Mittag-Leffler condition with respect to the sequence $(l(m))_{m \in \mathbb{N}}$ if and only if

$$\Lambda_m \cdots \Lambda_k(H_{k+1}) = \Lambda_m \cdots \Lambda_{l(m)-1}(H_{l(m)}) \qquad (*)$$

for any $k \geq l(m)$. By Lemma 1.1, the composition of diagonal maps is a diagonal map and we can conclude that the equality $(*)$ is equivalent to

$$\lambda_m^i \cdots \lambda_k^i(H_{k+1}) = \lambda_m^i \cdots \lambda_{l(m)-1}^i(H_{l(m)}) \text{ for each } i \in I.$$

That is, for each $i \in I$, $\mathcal{T}_i$ satisfies the Mittag-Leffler condition with respect to the sequence $(l(m))_{m \in \mathbb{N}}$.

Again Lemma 1.1 allows us to conclude that the equality $(*)$ is equivalent to

$$\Lambda_m' \cdots \Lambda_k'(H_{k+1}) = \Lambda_m' \cdots \Lambda_{l(m)-1}'(H_{l(m)}).$$

That is, the tower $\bigoplus_{i \in I} \mathcal{T}_i$ satisfies the Mittag-Leffler condition with respect to the sequence $(l(m))_{m \in \mathbb{N}}$. ■

In [13], [17] and [7], the Mittag-Leffler condition for a tower of right $R$-modules $\mathcal{T} = (H_n, \lambda_n)_{n \in \mathbb{N}}$ is interpreted in terms of vanishing of $\varprojlim^1$, the first derived functor of the inverse limit $\varprojlim$. Recall that $\varprojlim^1$ is defined by the exact sequence

$$(\sharp) \qquad 0 \to \varprojlim H_n \to \prod_{n \in \mathbb{N}} H_n \xrightarrow{\Delta_{\mathcal{T}}} \prod_{n \in \mathbb{N}} H_n \to \varprojlim^1 H_n \to 0$$

where $\Delta_{\mathcal{T}}(a_n)_{n \in \mathbb{N}} = (a_n - \lambda_n(a_{n+1}))_{n \in \mathbb{N}}$ for any $(a_n)_{n \in \mathbb{N}} \in \prod_{n \in \mathbb{N}} H_n$, see [21, 3.5].

Thus, $\varprojlim^1 H_n = 0$ if and only if $\Delta_{\mathcal{T}}$ is surjective.

It has been known for a long time that the Mittag-Leffler condition is a sufficient condition for the vanishing of $\varprojlim^1$ [21, Proposition 3.5.7]. In [7], Emmanouil gave a necessary condition for a tower to be Mittag-Leffler in terms of vanishing of $\varprojlim^1$. Our next result follows from [7, Corollary 6], however we present an alternative proof, more direct than Emmanouil's one.

**Theorem 1.3.** *Let $\mathcal{T} = (H_n, \lambda_n)_{n \in \mathbb{N}}$ be a tower of right $R$-modules. Then the following statements are equivalent.*

(1) $\varprojlim^1 H_n^{(I)} = 0$ *for any set $I$,*
(2) $\varprojlim^1 H_n^{(\mathbb{N})} = 0$,
(3) $\mathcal{T} = (H_n, \lambda_n)_{n \in \mathbb{N}}$ *satisfies the Mittag-Leffler condition.*

PROOF. Clearly (1) implies (2). Assume (3). Since $\mathcal{T}$ satisfies the Mittag-Leffler condition, so does $\mathcal{T}^{(I)}$ for any set $I$. Then, by [21, Proposition 3.5.7], it follows that $\varprojlim^1 H_n^{(I)} = 0$ for any set $I$ as claimed in (1).

Now we prove that (2) implies (3). We identify $\Delta_\mathcal{T}$ and $\Delta_{\mathcal{T}^{(\mathbb{N})}}$ with the matrices

$$\Delta_\mathcal{T} = \begin{pmatrix} 1 & -\lambda_1 & 0 & \ldots & 0 & \ldots \\ 0 & 1 & -\lambda_2 & & & \\ \vdots & & \ddots & \ddots & & \vdots \\ 0 & & & 1 & -\lambda_n & \\ \vdots & & & & \ddots & \ddots \end{pmatrix} \quad \Delta_{\mathcal{T}^{(\mathbb{N})}} = \begin{pmatrix} 1 & -\Lambda'_1 & 0 & \ldots & 0 & \ldots \\ 0 & 1 & -\Lambda'_2 & & & \\ \vdots & & \ddots & \ddots & & \vdots \\ 0 & & & 1 & -\Lambda'_n & \\ \vdots & & & & \ddots & \ddots \end{pmatrix}$$

where, for each $n \in \mathbb{N}$, $\Lambda'_n : H_n^{(\mathbb{N})} \to H_n^{(\mathbb{N})}$ is the diagonal map relative to the sequence constantly equal to $\lambda_n$. That is, if $(a_i)_{i \in \mathbb{N}} \in H_{n+1}^{(\mathbb{N})}$ then $\Lambda'_n((a_i)_{i \in \mathbb{N}}) = (\lambda_n(a_i))_{i \in \mathbb{N}}$.

**Step 1.** *Let $B = (b_{ij})_{i\,j \in \mathbb{N}}$ be a row finite matrix, with $b_{ij} \in H_i$, for any $i, j \in \mathbb{N}$. Then (2) implies that there exists a row finite matrix $A = (a_{ij})_{i\,j \in \mathbb{N}}$, with $a_{ij} \in H_i$, for any $i, j \in \mathbb{N}$, such that $\Delta_\mathcal{T} A = B$.*

To prove Step 1, let $B_i \in H_i^{(\mathbb{N})}$ be the $i$-th row of $B$, for each $i \in \mathbb{N}$.

Since $\varprojlim^1 H_n^{(\mathbb{N})} = 0$ if and only if $\Delta_{\mathcal{T}^{(\mathbb{N})}}$ is onto, there exists a sequence $(A_i)_{i \in \mathbb{N}}$, $A_i = (a_{ij})_{j \in \mathbb{N}} \in H_i^{(\mathbb{N})}$ for each $i \in \mathbb{N}$, such that

$$\Delta_{\mathcal{T}^{(\mathbb{N})}} \begin{pmatrix} A_1 \\ A_2 \\ \vdots \\ A_n \\ \vdots \end{pmatrix} = \begin{pmatrix} 1 & -\Lambda'_1 & 0 & \ldots & 0 & \ldots \\ 0 & 1 & -\Lambda'_2 & & & \\ \vdots & & \ddots & \ddots & & \vdots \\ 0 & & & 1 & -\Lambda'_n & \\ \vdots & & & & \ddots & \ddots \end{pmatrix} \begin{pmatrix} A_1 \\ A_2 \\ \vdots \\ A_n \\ \vdots \end{pmatrix} = \begin{pmatrix} B_1 \\ B_2 \\ \vdots \\ B_n \\ \vdots \end{pmatrix}.$$

Equivalently, $\Delta_\mathcal{T} A = B$ where $A = (a_{ij})_{i\,j \in \mathbb{N}}$. This shows our claim.

**Step 2.** *Let $A$ and $B$ be matrices like in Step 1 such that $\Delta_\mathcal{T} A = B$. Assume $b_{ij} = 0$ for $i \neq j$. Then there exists a sequence of natural numbers $(l(m))_{m \in \mathbb{N}}$, with $l(m) > m$ for every $m \in \mathbb{N}$, satisfying that*

$$\lambda_m \cdots \lambda_{k-1}(b_{k\,k}) \in \lambda_m \cdots \lambda_k(H_{k+1})$$

*for all $k \geq l(m)$.*

The proof follows from an argument that goes back to Bass [2], see also the paper by Azumaya [1]. For the computation in this type of situation see, for example, the proof of [4, Lemma 3.3].

**Step 3.** $\mathcal{T} = (H_n, \lambda_n)_{n \in \mathbb{N}}$ satisfies the Mittag-Leffler condition.

Assume by way of contradiction that there exists an integer $m$ for which the chain
$$H_m \supseteq \lambda_m(H_{m+1}) \supseteq \cdots \supseteq \lambda_m \lambda_{m+1} \cdots \lambda_{n-1}(H_n) \supseteq \cdots$$
is not stationary. This means that there exists an infinite set $N \subseteq \mathbb{N}$ such that, for any $n \in N$, there is an element $b_n \in H_n$ such that $\lambda_m \cdots \lambda_{n-1}(b_n) \notin \lambda_m \cdots \lambda_n(H_{n+1})$. Consider the matrix $B = (b_{ij})_{i\,j \in \mathbb{N}}$ such that $b_{nn} = b_n$, for $n \in N$, and $b_{ij} = 0$ otherwise. By Step 1 there exists $A$ such that $\Delta_\mathcal{T} A = B$. By Step 2, there exists an integer $l(m) > m$ such that for all $k \geq l(m)$, $\lambda_m \cdots \lambda_{k-1}(b_{k\,k}) \in \lambda_m \cdots \lambda_k(H_{k+1})$ which contradicts the choice of the infinite family $(b_n)_{n \in N}$. ∎

**Examples 1.4.** *Let $\mathcal{T} = (H_n, \lambda_n)_{n \in \mathbb{N}}$ be a tower of right modules over a ring $R$. If either the homomorphisms $\lambda_n$ are onto or $H_n$ are artinian modules, then $\mathcal{T}$ satisfies the Mittag-Leffler condition.*

*Another trivial example of Mittag-Leffler tower is given by $\mathcal{T}$-nilpotent sequences of maps $(\lambda_n)_{n \in \mathbb{N}}$. In this case $\Delta_\mathcal{T}$ is an isomorphism.*

## 2. Countable direct limits of modules

Let $R$ be a ring. We fix the following notation.

**Notation 2.1.** *Given a countable direct system*
$$C_1 \xrightarrow{f_1} C_2 \xrightarrow{f_2} C_3 \to \ldots \to C_n \xrightarrow{f_n} C_{n+1} \to \ldots$$
*of right $R$-modules, we consider the pure exact sequence*
$$0 \to \oplus_{n \in \mathbb{N}} C_n \xrightarrow{\phi} \oplus_{n \in \mathbb{N}} C_n \to \varinjlim C_n \to 0$$
*where $\phi \varepsilon_n = \varepsilon_n - \varepsilon_{n+1} f_n$ and $\varepsilon_n \colon C_n \to \oplus_{n \in \mathbb{N}} C_n$ denotes the canonical morphism for every $n \in \mathbb{N}$.*

Let $M$ be a right $R$-module. Applying the functor $\mathrm{Hom}_R(-, M)$ to the setting of Notation 2.1 we obtain a tower of modules over $\mathrm{End}_R(M)$.

**Lemma 2.2.** *Consider the setting of Notation 2.1, and let $M$ be a right $R$-module with $S = \mathrm{End}_R(M)$. Then*
$$\mathcal{T} = (\mathrm{Hom}_R(C_n, M), \mathrm{Hom}_R(f_n, M))_{n \in \mathbb{N}}$$
*is a tower of left $S$-modules. Moreover, the following diagram is commutative*

$$\begin{array}{ccc} \mathrm{Hom}_R(\bigoplus_{n \in \mathbb{N}} C_n, M) & \xrightarrow{\mathrm{Hom}_R(\phi, M)} & \mathrm{Hom}_R(\bigoplus_{n \in \mathbb{N}} C_n, M) \\ \downarrow \cong & & \cong \downarrow \\ \prod_{n \in \mathbb{N}} \mathrm{Hom}_R(C_n, M) & \xrightarrow{\Delta_\mathcal{T}} & \prod_{n \in \mathbb{N}} \mathrm{Hom}_R(C_n, M) \end{array}$$

*where $\Delta_\mathcal{T}$ is the homomorphism in the sequence $(\sharp)$.*

We now translate the results of § 1 to this particular kind of towers. Recall that a right $R$-module $C$ is said to be *small* if $\mathrm{Hom}_R(C, \oplus_{i \in I} M_i) \cong \oplus_{i \in I} \mathrm{Hom}_R(C, M_i)$ for every family of right $R$-modules $\{M_i \mid i \in I\}$. Finitely generated modules are examples of small modules.

**Corollary 2.3.** *Let $R$ be a ring and $M$ a right $R$-module with $\mathrm{End}_R(M) = S$. Assume that the modules in Notation 2.1 are small. Then the following statements are equivalent:*

(1) $\mathrm{Hom}_R(\phi, M^{(I)})$ *is surjective for any set $I$,*
(2) $\mathrm{Hom}_R(\phi, M^{(\mathbb{N})})$ *is surjective,*
(3) $\varprojlim^1 \mathrm{Hom}_R(C_n, M)^{(\mathbb{N})} = 0$,
(4) *The tower $(\mathrm{Hom}_R(C_n, M), \mathrm{Hom}_R(f_n, M))_{n \in \mathbb{N}}$ satisfies the Mittag-Leffler condition,*
(5) *For every diagonal map $\gamma : \bigoplus_{n \in \mathbb{N}} C_n \to M^{(\mathbb{N})}$ there is a homomorphism $\psi : \bigoplus_{n \in \mathbb{N}} C_n \to M^{(\mathbb{N})}$ such that $\psi \phi = \gamma$.*

PROOF. Let $I$ be a set, and let $\mathcal{T} = (\mathrm{Hom}_R(C_n, M), \mathrm{Hom}_R(f_n, M))_{n \in \mathbb{N}}$. Since all modules $C_n$ are small, the tower $(\mathrm{Hom}_R(C_n, M^{(I)}), \mathrm{Hom}_R(f_n, M^{(I)}))_{n \in \mathbb{N}}$ is naturally isomorphic to the tower $\mathcal{T}^{(I)}$. Then we know from Lemma 2.2 that $\mathrm{Hom}_R(\phi, M^{(I)})$ coincides up to natural isomorphism with $\Delta_{\mathcal{T}^{(I)}}$. So the equivalence of (1), (2), (3), (4) follows from Theorem 1.3. Clearly, condition (2) implies (5); the implication (5) $\Rightarrow$ (4) is exactly the implication (i) $\Rightarrow$ (iii) in [4, Theorem 3.7] . ∎

**Examples 2.4.** We assume the setting of Notation 2.1 where the modules $(C_n)_{n \in \mathbb{N}}$ are small.
(1) If the sequence

$$0 \to \oplus_{n \in \mathbb{N}} C_n \xrightarrow{\phi} \oplus_{n \in \mathbb{N}} C_n \to \varinjlim C_n \to 0 \qquad (*)$$

splits, then $\mathrm{Hom}_R(\phi, M^{(I)})$ is surjective for any module $M$. Hence, by Corollary 2.3, the tower $(\mathrm{Hom}_R(C_n, M), \mathrm{Hom}_R(f_n, M))_{n \in \mathbb{N}}$ satisfies the Mittag-Leffler condition for any module $M$.
The converse is also true. Assume the tower $(\mathrm{Hom}_R(C_n, M), \mathrm{Hom}_R(f_n, M))_{n \in \mathbb{N}}$ satisfies the Mittag-Leffler condition for any module $M$. Taking $M = \oplus_{n \in \mathbb{N}} C_n$, we deduce from Corollary 2.3 that $\mathrm{Hom}_R(\phi, (\oplus_{n \in \mathbb{N}} C_n)^{(\mathbb{N})})$ is surjective and, hence, $\mathrm{Hom}_R(\phi, \oplus_{n \in \mathbb{N}} C_n)$ is surjective. Thus, there exists $\psi : \oplus_{n \in \mathbb{N}} C_n \to \oplus_{n \in \mathbb{N}} C_n$ such that $\psi \circ \phi = \mathrm{Id}$. This shows that $(*)$ splits.
(2) If $M$ is a $\Sigma$-pure-injective module, then $\mathrm{Hom}_R(\phi, M^{(I)})$ is surjective for any set $I$. It follows from Corollary 2.3 that the tower $(\mathrm{Hom}_R(C_n, M), \mathrm{Hom}_R(f_n, M))_{n \in \mathbb{N}}$ satisfies the Mittag-Leffler condition.
(3) If $M$ is injective and the maps $f_n$ are monomorphisms, then $\mathrm{Hom}_R(f_n, M)$ is surjective for any $n \in \mathbb{N}$, hence $(\mathrm{Hom}_R(C_n, M), \mathrm{Hom}_R(f_n, M))_{n \in \mathbb{N}}$ clearly satisfies the Mittag-Leffler condition.
(4) If $M$ is a module such that

$$\mathrm{Ext}^1_R(\varinjlim C_n, M^{(\mathbb{N})}) = 0$$

then $\mathrm{Hom}_R(\phi, M^{(\mathbb{N})})$ is surjective. Hence, by Corollary 2.3, $(\mathrm{Hom}_R(C_n, M), \mathrm{Hom}_R(f_n, M))_{n \in \mathbb{N}}$ satisfies the Mittag-Leffler condition. If, moreover, $\mathrm{Ext}^1_R(C_n, M^{(\mathbb{N})}) = 0$ for any $n \in \mathbb{N}$, then the converse is also true (cf. [4, Theorem 5.1]).

We are interested in further developing Example 2.4(4) in the case where each $C_n$ is a finitely generated free module.

**Proposition 2.5.** *Assume that the modules in Notation 2.1 are finitely generated and free. Then $\mathrm{Ext}^1_R(\varinjlim C_n, M^{(\mathbb{N})}) = 0$ if and only if $(\mathrm{Hom}_R(C_n, M), \mathrm{Hom}_R(f_n, M))_{n \in \mathbb{N}}$ satisfies the Mittag-Leffler condition.*

*In particular, $\varinjlim C_n$ is projective if and only if $(\mathrm{Hom}_R(C_n, R), \mathrm{Hom}_R(f_n, R))_{n \in \mathbb{N}}$ satisfies the Mittag-Leffler condition.*

PROOF. As $\oplus_{n \in \mathbb{N}} C_n = F$ is a free module, applying the functor $\mathrm{Hom}_R(-, M^{(\mathbb{N})})$ to the exact sequence
$$0 \to \oplus_{n \in \mathbb{N}} C_n = F \xrightarrow{\phi} \oplus_{n \in \mathbb{N}} C_n = F \to \varinjlim C_n \to 0 \qquad (*)$$
we obtain the exact sequence
$$0 \to \mathrm{Hom}_R(\varinjlim C_n, M^{(\mathbb{N})}) \to$$
$$\to \mathrm{Hom}_R(F, M^{(\mathbb{N})}) \xrightarrow{\mathrm{Hom}_R(\phi, M^{(\mathbb{N})})} \mathrm{Hom}_R(F, M^{(\mathbb{N})}) \to \mathrm{Ext}^1_R(\varinjlim C_n, M^{(\mathbb{N})}) \to 0.$$
This shows that $\mathrm{Ext}^1_R(\varinjlim C_n, M^{(\mathbb{N})}) = 0$ if and only if $\mathrm{Hom}_R(\phi, M^{(\mathbb{N})})$ is surjective. By Corollary 2.3, this is equivalent to the fact that $(\mathrm{Hom}_R(C_n, M), \mathrm{Hom}_R(f_n, M))_{n \in \mathbb{N}}$ satisfies the Mittag-Leffler condition.

Assume now that $(\mathrm{Hom}_R(C_n, R), \mathrm{Hom}_R(f_n, R))_{n \in \mathbb{N}}$ satisfies the Mittag-Leffler condition. Then, by Corollary 2.3, $\mathrm{Hom}_R(\phi, F)$ is surjective. Hence, there exists $\psi \colon F \to F$ such that $\psi \circ \phi = \mathrm{Id}$. This implies that the sequence $(*)$ splits, and we conclude that $\varinjlim C_n$ is projective. ∎

We now collect some closure properties of the class of all modules $M$ that turn the direct system of Notation 2.1 into a Mittag-Leffler tower $(\mathrm{Hom}_R(C_n, M), \mathrm{Hom}_R(f_n, M))_{n \in \mathbb{N}}$.

**Corollary 2.6.** *Assume that the modules in Notation 2.1 are small. For a given set $I$, let $\{M_i\}_{i \in I}$ be a family of right $R$ modules. Then the tower*
$$\left( \mathrm{Hom}_R(C_n, \bigoplus_{i \in I} M_i), \mathrm{Hom}_R(f_n, \bigoplus_{i \in I} M_i) \right)_{n \in \mathbb{N}}$$
*satisfies the Mittag-Leffler condition if and only if so does the tower*
$$\left( \mathrm{Hom}_R(C_n, \prod_{i \in I} M_i), \mathrm{Hom}_R(f_n, \prod_{i \in I} M_i) \right)_{n \in \mathbb{N}}.$$

PROOF. The claim follows from Proposition 1.2 since there are isomorphisms of towers $(\mathrm{Hom}_R(C_n, \prod_{i \in I} M_i), \mathrm{Hom}_R(f_n, \prod_{i \in I} M_i))_{n \in \mathbb{N}} \cong \prod_{i \in I} (\mathrm{Hom}_R(C_n, M_i), \mathrm{Hom}_R(f_n, M_i))_{n \in \mathbb{N}}$ and $(\mathrm{Hom}_R(C_n, \oplus_{i \in I} M_i), \mathrm{Hom}_R(f_n, \oplus_{i \in I} M_i))_{n \in \mathbb{N}} \cong \bigoplus_{i \in I} (\mathrm{Hom}_R(C_n, M_i), \mathrm{Hom}_R(f_n, M_i))_{n \in \mathbb{N}}$. ∎

The next closure property relies on [4, Lemma 4.1], a result which we recall in 2.8. Hereby, we present a more elegant proof which was suggested to us by the referee and which uses the following "homotopy lemma", cf. [15, Lemma B1, Appendix B].

**Lemma 2.7.** *Consider the commutative diagram of right modules and module homomorphisms*

$$\begin{array}{ccccccc} C & \xrightarrow{f} & C' & \xrightarrow{\pi} & C'' & \longrightarrow & 0 \\ {\scriptstyle h}\downarrow & & {\scriptstyle k}\downarrow & & {\scriptstyle \ell}\downarrow & & \\ 0 & \longrightarrow & N & \xrightarrow{\varepsilon} & M & \xrightarrow{g} & L \end{array}$$

*and assume it has exact rows.*

*Then there exists $q\colon C'' \to M$ such that $gq = \ell$ if and only if there exists $p\colon C' \to N$ such that $pf = h$.*

PROOF. Let $q\colon C'' \to M$ be such that $gq = \ell$. As $\ell\pi = gk$, $g(k - q\pi) = 0$. Since $\varepsilon\colon N \to M$ is the kernel of $g$, there exists $p\colon C' \to N$ such that $\varepsilon p = k - q\pi$. Then $\varepsilon pf = (k - q\pi)f = kf = \varepsilon h$. As $\varepsilon$ is a monomorphism, we deduce that $pf = h$.

The converse follows by a dual argument. ∎

**Lemma 2.8.** *([4, Lemma 4.1]) Let $C$ and $C'$ be finitely generated right $R$-modules such that $C'$ is finitely presented, and let $f\colon C \to C'$ be a module homomorphism. If $M$ is a right $R$-module with a pure submodule $N$ then*

$$\operatorname{Hom}_R(C', M)f \cap \operatorname{Hom}_R(C, N) = \operatorname{Hom}_R(C', N)f.$$

PROOF. Denote by $\varepsilon\colon N \to M$ the inclusion. Let $\pi\colon C' \to \operatorname{Coker} f$ and $g\colon M \to \operatorname{Coker} \varepsilon$ denote the cokernel of $f$ and the cokernel of $\varepsilon$, respectively.

It is clear that $\operatorname{Hom}_R(C', N)f \subseteq \operatorname{Hom}_R(C', M)f \cap \operatorname{Hom}_R(C, N)$. Assume $k \in \operatorname{Hom}_R(C', M)$ is such that $kf(C) \subseteq N$, and set $h = kf\colon C \to N$. As $gkf = g\varepsilon h = 0$, there exists $\ell\colon \operatorname{Coker} f \to \operatorname{Coker} \varepsilon$ such that $\ell\pi = gk$. Therefore we have the commutative diagram with exact rows

$$\begin{array}{ccccccc} C & \xrightarrow{f} & C' & \xrightarrow{\pi} & \operatorname{Coker} f & \longrightarrow & 0 \\ {\scriptstyle h}\downarrow & & {\scriptstyle k}\downarrow & & {\scriptstyle \ell}\downarrow & & \\ 0 & \longrightarrow & N & \xrightarrow{\varepsilon} & M & \xrightarrow{g} & \operatorname{Coker} \varepsilon & \longrightarrow & 0 \end{array}$$

As $C$ is finitely generated and $C'$ is finitely presented, $\operatorname{Coker} f$ is a finitely presented module. Then, since the sequence

$$0 \to N \xrightarrow{\varepsilon} M \xrightarrow{g} \operatorname{Coker} \varepsilon \to 0$$

is pure-exact, the sequence

$$\operatorname{Hom}_R(\operatorname{Coker} f, M) \xrightarrow{\operatorname{Hom}_R(\operatorname{Coker} f, g)} \operatorname{Hom}_R(\operatorname{Coker} f, \operatorname{Coker} \varepsilon) \to 0$$

is exact. In particular, there exists $q \in \operatorname{Hom}_R(\operatorname{Coker} f, M)$ such that $gq = \ell$. By Lemma 2.7, there exists $p\colon C' \to N$ such that $h = pf$. This finishes the proof of the statement. ∎

**Proposition 2.9.** *Assume that the modules $C_n$ in the Notation 2.1 are finitely presented. Let $M$ be a right $R$-module such that the tower $(\operatorname{Hom}_R(C_n, M), \operatorname{Hom}_R(f_n, M))_{n \in \mathbb{N}}$ satisfies the Mittag-Leffler condition. Then, for every pure submodule $N$ of $M$ the tower $(\operatorname{Hom}_R(C_n, N), \operatorname{Hom}_R(f_n, N))_{n \in \mathbb{N}}$ satisfies the Mittag-Leffler condition.*

PROOF. Let $N$ be a pure submodule of $M$. By Lemma 2.8, for $m, n \in \mathbb{N}$ and any map $f \in \operatorname{Hom}_R(C_m, C_{m+n})$

$$\operatorname{Hom}_R(C_{m+n}, M)f \cap \operatorname{Hom}_R(C_m, N) = \operatorname{Hom}_R(C_{m+n}, N)f.$$

So, if the chain of subgroups of $\operatorname{Hom}_R(C_m, M)$

$$\operatorname{Hom}_R(C_{m+1}, M)f_m \supseteq \operatorname{Hom}_R(C_{m+2}, M)f_{m+1}f_m \supseteq \ldots$$
$$\cdots \supseteq \operatorname{Hom}_R(C_{m+n}, M)f_{m+n-1}f_{m+n-2}\cdots f_m \supseteq \ldots$$

is stationary, also the corresponding chain of subgroups of $\operatorname{Hom}_R(C_m, N)$ is stationary. This proves the claim. ∎

## 3. BAER MODULES ARE PROJECTIVE

We will now apply the previous results to the Baer splitting problem. As observed in the introduction, we will only have to consider Baer modules that are countably presented and flat. For such modules we have the following result which is essentially well known; the idea goes back to Jensen's proof of the fact that countably presented flat modules have projective dimension at most one [14, Lemma 2]. We sketch the argument for sake of completeness.

**Proposition 3.1.** *Let $B$ be a countably presented flat right module over a ring $R$. Then there is a countable direct system*

$$F_1 \xrightarrow{f_1} F_2 \xrightarrow{f_2} F_3 \to \ldots \to F_n \xrightarrow{f_n} F_{n+1} \to \ldots$$

*where $F_n$ are finitely generated free modules such that, following Notation 2.1, $B$ fits in the exact sequence*

$$0 \to \oplus_{n \in \mathbb{N}} F_n \xrightarrow{\phi} \oplus_{n \in \mathbb{N}} F_n \to \varinjlim F_n = B \to 0 \qquad (*).$$

PROOF. By hypothesis, the module $B$ has a presentation of the form

$$0 \to G \xrightarrow{f} F \to B \to 0,$$

where $F$ is a countably generated free right module and $G$ is a countably generated right module. Let $\{g_n\}_{n \in \mathbb{N}}$ be a set of generators of $G$, and let $\{e_n\}_{n \in \mathbb{N}}$ be a basis of $F$. For a set $A \subseteq \mathbb{N}$ we denote by $F_A = \sum_{n \in A} e_n R$.

Following Jensen [14, proof of Lemma 2], we can find an ascending chain $(A_n)_{n \in \mathbb{N}}$ of finite subsets of $\mathbb{N}$ such that
(1) $\bigcup_{n \in \mathbb{N}} A_n = \mathbb{N}$,
(2) for any $i \leq n$, $f(g_i) \in F_{A_n}$,
and, as $f$ is a pure monomorphism, also
(3) the induced map $\bar{f} \colon \sum_{i \leq n} g_i R \to F_{A_n}$ splits.

Condition (3) implies that $P_n = F_{A_n} / \sum_{i \leq n} f(g_i) R$ is finitely generated and projective.

For each $n \in \mathbb{N}$, let $f'_n \colon P_n \to P_{n+1}$ be the map induced by $f$ and choose a finitely generated projective module $Q_n$ such that $P_n \oplus Q_n = F_n$ is finitely generated and free. Finally, for each $n \in \mathbb{N}$, let $f_n \colon F_n \to F_{n+1}$ be the homomorphism defined by $f_n(p, q) = (f'_n(p), 0)$, for $(p, q) \in P_n \oplus Q_n = F_n$.

Then we obtain a direct system
$$F_1 \xrightarrow{f_1} F_2 \xrightarrow{f_2} F_3 \to \ldots \to F_n \xrightarrow{f_n} F_{n+1} \to \ldots$$
and it is easy to check that its direct limit is $B$. ∎

We finally specialize to commutative domains. We start with some preliminary results.

**Lemma 3.2.** *Let $R$ be a commutative domain that is not a field. Let $G$ be a finitely generated $R$-module. Then $\bigcap_{0 \neq r \in R} rG = 0$.*

PROOF. We first assume that $G$ is finitely generated and torsion-free. In this case $\bigcap_{0 \neq r \in R} rG$ coincides with the divisible submodule $d(G)$ of $G$ which is torsion-free and divisible, hence isomorphic to a direct sum of copies of $Q$. If $d(G) \neq 0$, $G$ contains a summand isomorphic to $Q$ which is impossible since $G$ is finitely generated.

If $G$ is not torsion-free, consider the exact sequence:
$$0 \to t(G) \to G \to G/t(G) \to 0$$
where $t(G)$ denotes the torsion submodule of $G$. Let $x \in \bigcap_{0 \neq r \in R} rG$; then $x + t(G) \in \bigcap_{0 \neq r \in R} (rG + t(G))/t(G) = \bigcap_{0 \neq r \in R} r(G/t(G))$ which is zero by the first part of the proof. So $x \in t(G)$. Moreover, since $rG \cap t(G) = rt(G)$ for every $r \in R$, we have that $x \in \bigcap_{0 \neq r \in R} rt(G)$. By induction on the number of generators, it is not difficult to show that the torsion submodule of a finitely generated module has nonzero annihilator (cf. [18, Lemma 7.1]). Thus $\bigcap_{0 \neq r \in R} rt(G) = 0$; hence $x = 0$. ∎

**Lemma 3.3.** *Let $R$ be a commutative domain that is not a field. Then the map*
$$\mu \colon R \to \prod_{0 \neq r \in R} R/rR$$
*defined by $\mu(x) = (x + rR)$, for any $x \in R$, is a pure embedding.*

PROOF. By Lemma 3.2, $\mu$ is a monomorphism. We show that $\mu$ is a pure monomorphism. We prove that for every finitely presented $R$-module $G$, the homomorphism:
$$\mu \otimes_R 1_G \colon R \otimes_R G \to \prod_{0 \neq r \in R} R/rR \otimes_R G$$
is a monomorphism. Since $G$ is finitely presented, $- \otimes_R G$ commutes with direct products, so we are lead to show that the homomorphism
$$\nu \colon G \to \prod_{0 \neq r \in R} G/rG$$

is a monomorphism. As $\operatorname{Ker} \nu = \bigcap_{0 \neq r \in R} rG$, we conclude by Lemma 3.2. ∎

We are now in a position to prove our main result.

**Theorem 3.4.** *Let $B$ be a Baer module over a commutative domain $R$. Then $B$ is projective.*

PROOF. As proved in [8] a module $B$ is a Baer module if and only if there exists a well ordered continuous ascending chain $(B_\alpha \mid \alpha < \lambda)$ of submodules such that all the factors $B_{\alpha+1}/B_\alpha$ are countably generated Baer modules. Thus, by [6, XII Lemma 1.5], to decide whether Baer modules are projective it is enough to consider countably generated Baer modules.

Let $B$ be a countably generated Baer module. By Kaplansky's results [16], $B$ is flat and of projective dimension at most one. Recall that over a commutative domain, countably generated flat modules have projective dimension at most one if and only if they are countably presented. Hence $B$ is a flat countably presented module, and we can fix a direct system of finitely generated free modules

$$F_1 \xrightarrow{f_1} F_2 \xrightarrow{f_2} F_3 \to \ldots \to F_n \xrightarrow{f_n} F_{n+1} \to \ldots$$

as given by Proposition 3.1 such that $\varinjlim F_n \cong B$.

Let $M = \oplus_{0 \neq r \in R} R/rR$. As $M$ and hence $M^{(\mathbb{N})}$ are torsion modules, $\operatorname{Ext}_R^1(B, M^{(\mathbb{N})}) = 0$. By Proposition 2.5, the tower

$$(\operatorname{Hom}_R(F_n, M), \operatorname{Hom}_R(f_n, M))_{n \in \mathbb{N}}$$

satisfies the Mittag-Leffler condition. By Corollary 2.6, the tower

$$\left( \operatorname{Hom}_R(F_n, \prod_{0 \neq r \in R} R/rR), \operatorname{Hom}_R(f_n, \prod_{0 \neq r \in R} R/rR) \right)_{n \in \mathbb{N}}$$

also satisfies the Mittag-Leffler condition.

By Lemma 3.3, $R$ is a pure submodule of $\prod_{0 \neq r \in R} R/rR$. Applying Proposition 2.9, we infer that the tower $(\operatorname{Hom}_R(F_n, R), \operatorname{Hom}_R(f_n, R))_{n \in \mathbb{N}}$ satisfies the Mittag-Leffler condition. Hence, by Proposition 2.5, we conclude that $B$ is projective. ∎

To prove that a countably generated Baer module $B$ is projective we have only used that $\operatorname{Ext}_R^1(B, (\oplus_{0 \neq r \in R} R/rR)^{(\mathbb{N})}) = 0$, but we note that this is in fact equivalent to the statement that $B$ is a Baer module (cf. [9, Proposition 8.14] or [5]). In general, we obtain as a consequence

**Corollary 3.5.** *Let $R$ be a commutative domain, and let $\kappa$ be an infinite cardinal. Let $B$ be an $R$-module that can be generated by a set of cardinality at most $\kappa$. Then $B$ is projective if and only if*

$$\operatorname{Ext}_R^1(B, (\oplus_{0 \neq r \in R} R/rR)^{(\kappa)}) = 0$$

PROOF. By [9, Proposition 8.14], the statement is equivalent to say that $B$ is a Baer module. The conclusion follows from Theorem 3.4. ∎

Dipartimento di Informatica e Comunicazione, Università degli Studi dell'Insubria, Via Mazzini 5, I - 21100 Varese, Italy

*E-mail address*: lidia.angeleri@uninsubria.it

Dipartimento di Matematica Pura e Applicata, Università di Padova, Via Belzoni 7, I-35131 Padova, Italy

*E-mail address*: bazzoni@math.unipd.it

Departament de Matemàtiques, Universitat Autònoma de Barcelona, E-08193 Bellaterra (Barcelona), Spain

*E-mail address*: dolors@mat.uab.es